\pdfoutput=1
\date{} % %
\title{Planar Ising magnetization field II. Properties of the critical and near-critical scaling limits}
\author{Federico Camia \and Christophe Garban \and Charles M. Newman}

\documentclass[11pt,preprint]{article}

\usepackage[vmargin={3cm, 3.5cm},hmargin=2.5cm]{geometry}

\usepackage[pdfpagelabels]{hyperref}

\usepackage{latexsym} 
\usepackage{amsmath}
\usepackage{amssymb}
\usepackage{amsthm}
\usepackage{mathrsfs}
\usepackage{amsfonts}
\usepackage{graphicx}
\usepackage{url}

\usepackage{color}
\usepackage{bbm,verbatim}

\usepackage{minitoc} %% POUR LES TABLES MATIERES DANS LES CHAPITRES.

\usepackage{multirow} %pour les tableaux
\usepackage{array}
\newcolumntype{M}[1]{>{\centering}m{#1}}

\numberwithin{equation}{section}
\numberwithin{figure}{section}
\newtheorem{theorem}{Theorem}
\numberwithin{theorem}{section}
\newtheorem{lemma}[theorem]{Lemma}

\newtheorem{proposition}[theorem]{Proposition}

\theoremstyle{remark}\newtheorem{remark}[theorem]{Remark}
\def\eqref#1{(\ref{#1})}
\let\qqed=\qed

\def\nn{\nonumber}

\def\QED{\qqed\medskip}
\let\qed=\QED

\newcommand{\R}{\mathbb{R}}

\newcommand{\C}{\mathbb{C}}
\newcommand{\Z}{\mathbb{Z}}
\newcommand{\N}{\mathbb{N}}

\def\SLEkk#1/{$\mathrm{SLE}(#1)$}
\def\SLEr#1/{$\mathrm{SLE(\kappa;#1)}$}
\def\SLEkr#1;#2/{$\mathrm{SLE(#1;#2)}$}
\def\SLEk/{\SLEkk{\kappa}/}
\def\SLEtwo/{\SLEkk2/}
\def\SLE/{$\mathrm{SLE}$}
\def\SLEab/{\SLEkr 4; {a/\hco-1}, {b/\hco-1}/}

\def \eps {\epsilon}
\def \P {\mathbb{P}}
\def \E {{\mathbb E}}

\def\md{\mid}
\def\Bb#1#2{{\def\md{\bigm| }#1\bigl[#2\bigr]}}

\def\Pb{\Bb\P}
\def\Eb{\Bb\E}

\def \p {{\partial}}

\def\FK#1#2#3{{\def\md{\bigm| } \P_{#1}^{\,#2}  \bigl[  #3 \bigr]}}
\def\EFK#1#2#3{{\def\md{\bigm| } \E_{#1}^{\,#2}  \bigl[  #3 \bigr]}}

\def\fk{\mathrm{FK}}
\def\free{\mathrm{free}}
\def\wired{\mathrm{wired}}

\def \p {{\partial}}

\def\1{1\hspace{-2.55 mm}{1}}

\def\noopsort#1{}

\usepackage{mathrsfs}

\def\ni{\noindent}
\def\bi{\begin{itemize}}
\def\ei{\end{itemize}}

\def\Sob{\mathcal{H}} %Sobolev Hilbert space
\def\magn{m}

\def\<#1{\langle #1\rangle}
\def\g{\mathbf{g}}

\def\F{\mathcal{F}}

\begin{document}
\maketitle

\begin{abstract}

In \cite{CGN1}, we proved that the renormalized critical
Ising magnetization fields $\Phi^a:= a^{15/8} \sum_{x\in a\, \Z^2} \sigma_x \, \delta_x$ converge as $a\to 0$ to a random distribution that we denoted by $\Phi^\infty$. 
The purpose of this paper is to establish some fundamental properties satisfied by this 
$\Phi^\infty$ and the near-critical fields $\Phi^{\infty,h}$. 
More precisely, we obtain the following results.
\bi
\item[(i)] If $A\subset \C$ is a smooth bounded domain and if $m=m_A := \<{\Phi^\infty, 1_A}$ denotes the limiting rescaled magnetization in $A$, then there is a constant $c=c_A>0$ such that  
\begin{equation*}
\log \Pb{m > x} \underset{x\to \infty}{\sim}  -c \; x^{16}\,.
\end{equation*}
In particular, this provides an alternative proof that the field $\Phi^\infty$ is non-Gaussian (another proof of this fact would use the $n$-point correlation functions established in \cite{CHI} which do not satisfy Wick's formula). 
\item[(ii)] The random variable $m=m_A$ has a smooth {\it density} and one has more precisely the following bound on its Fourier transform: $|\Eb{e^{i\,t m}} |\le e^{- \tilde{c}\, |t|^{16/15}}$.
\item[(iii)] There exists a one-parameter family $\Phi^{\infty,h}$ of near-critical scaling limits for the
magnetization field in the plane with vanishingly small external magnetic field. 
\ei

\end{abstract}

%\tableofcontents

\section{Introduction}

\subsection{Overview}

In \cite{CGN1}, we considered the scaling limit of the appropriately renormalized magnetization field of a critical Ising model (i.e., at $\beta=\beta_c$) on the lattice $a\, \Z^2$ where the mesh $a$ shrinks to zero. The natural object to consider is the following field:
\begin{align}\label{}
\Phi^a:= \sum_{x\in a\, \Z^2}\, a^{15/8}\, \sigma_x\, \delta_x\,, 
\end{align}
where $\{\sigma_x\}_{x\in a\, \Z^2}$ is the realization of a critical Ising model on $a\, \Z^2$. Note 
that the renormalization of 
$a^{15/8}$ assumes Wu's Theorem of \cite{Wu}. See the introduction in \cite{CGN1} for a 
discussion of this. The following theorem is proved in \cite{CGN1}.

\begin{theorem}[\cite{CGN1}]\label{}
As the mesh $a\searrow 0$, the random field $\Phi^a$ converges in law to a limiting random field $\Phi^\infty$ under the topology of the Sobolev space $\Sob^{-3}(\C)$. See Theorem 1.2 and Appendix A in \cite{CGN1} for more details. 
\end{theorem}

In the case of a bounded smooth simply connected domain $\Omega$ equipped with $+,-$ or $\free$ boundary conditions along $\p \Omega$, one also obtains a limiting magnetization field $\Phi^\infty_{\Omega}$ whose law depends on the choice of the prescribed boundary conditions $\xi \in \{+,-,\free\}$. See Theorem 1.3 in \cite{CGN1}.

Two proofs of these results are provided in \cite{CGN1}: the first one relies on the recent breakthrough results from \cite{CHI} on the $n$-point correlation functions of the critical Ising model. The second proof is more conditional and relies for example on the on-going work \cite{IsingSLE}. See Section 2 in \cite{CGN1}.  
\medskip

The purpose of \cite{CGN1} was to identify a limit in law for these magnetization fields (i.e. $\Phi^\infty$ and $\Phi^\infty_\Omega$). Beyond the conformal covariance nature of these fields 
(Theorem 1.8 in \cite{CGN1}), we did not investigate the fine properties of these fields. This is what we wish to address in this paper:
\begin{enumerate}
\item  To start with, we will focus on the {\bf tail behavior} of the field $\Phi^\infty$ (and its bounded domain analog). For any bounded smooth domain $A$, we obtain a precise tail estimate for the block magnetization $m=m_A=\<{\Phi^\infty, 1_A}$ of $e^{-c\, x^{16}}$ where the constant $c>0$ does not depend on the prescribed boundary conditions along $\p A$. See Theorem \ref{th.tail}.

\item Then, we investigate whether the random variable $m_A$ defined above has a density function or not and if so what is its regularity. We answer this question by studying the tail of its characteristic function. Namely we prove that $|\Eb{e^{i\,t m}} |\le e^{- \tilde{c}\, |t|^{16/15}}$. See Theorem \ref{th.analyt}. 

\item Finally, we address a question of a different flavor: we prove in Theorem \ref{th.NC} that the magnetization field $\Phi^{a,h}$ for the near-critical Ising model with external field $h\, a^{15/8}$ has a scaling limit denoted by $\Phi^{\infty,h}$. 
%This limit is massive in that it has exponentially decaying correlations
%\margin{added parenthetical}
%(see a precise definition just before Theorem~\ref{th.NC} below).
%\margin{F: new sentence}
\end{enumerate}

\subsection{Main statements}

In Section \ref{s.tail} we will prove the following result.

\begin{theorem}\label{th.tail}
There exists a universal constant $c>0$ such that for any prescribed boundary conditions $\xi\in \{+,-,\free\}$ around the square $[0,1]^2$, 
the (continuum) magnetization $m=m^\xi= \Phi^\xi([0,1]^2)$ in $[0,1]^2$ satisfies as $x\to \infty$:

\[
\log \Pb{m> x} \sim -c\, x^{16}\,.
\]
This result extends to the case of the plane field $\Phi^\infty$ tested against a bounded smooth domain $A$, i.e.,
$\Phi^{\infty} (1_A)$ (in which case the constant $c$ will depend on $A$), or to the case of the limiting field $\Phi^\infty_\Omega$ for a bounded smooth domain $\Omega$ tested against a smooth sub-domain $A\subset \Omega$. 
\end{theorem}

%\begin{theorem}\label{th.tail}
%There exists a universal constant $c>0$ such that for any prescribed boundary conditions $\xi$ around the square $[0,1]^2$, 
%the (continuum) magnetization $m=m^\xi= \Phi^\xi([0,1]^2)$ in $[0,1]^2$ satisfies as $x\to \infty$:
%
%\[
%\log \Pb{m> x} \sim -c\, x^{16}\,.
%\]
%\end{theorem}

\ni
In Section \ref{s.analyt}, we will prove:

\begin{theorem}\label{th.analyt}
Let us consider the scaling limit $m=m^\xi$ of the magnetization in the square $[0,1]^2$ with  prescribed boundary conditions $\xi\in \{+,-,\free\}$.
% with boundary conditions $\xi$ around $\p [0,1]^24 in $\{+,-,\free\}$. 
There is a constant $\tilde{c}>0$ such that for all $t\in \R$ one has 
\[
|\EFK{}{\xi}{e^{i\, t\, m}}| \le e^{- \tilde{c} \, |t|^{\frac {16}{15}}}\,.
\]
%where the boundary conditions are either $+,-$ or $\free$. 
In particular, the density function $f=f^\xi$ of the random variable $m=m^\xi$ can be extended to an  entire function on the whole complex plane $\C$ \footnote{See for example Theorem IX.13 in \cite{Simon}}. 

As in Theorem \ref{th.tail}, the result extends to the 
whole-plane field $\Phi^\infty$ tested against domains $A$ as well as to the fields $\Phi^\infty_\Omega$ for smooth bounded domains $\Omega$. 
\end{theorem}

As we shall see later in Remark \ref{r.xi}, this theorem should also easily extend to more general boundary conditions $\xi$ such as finite combinations of $+,-,\free$-arcs. In this case, the constant $c=c([0,1]^2)>0$ would still be independent of the boundary condition $\xi$.    

%\begin{theorem}\label{th.analyt}
%There is a constant $c>0$ such that for all $t\in \R$ and for any boundary condition $\xi$, one has 
%\[
%|\EFK{}{\xi}{e^{i\, t\, m}}| \le e^{- c \, |t|^{\frac {16}{15}}}\,.
%\]
%In particular, the density function $f=f_\Omega$ is an {\bf entire} function on the whole plane $\C$.
%\end{theorem}

\bigskip

\ni
Finally, in Section \ref{s.NC} we will prove the following theorem concerning the near critical
(as $h\to 0$) scaling limit. Two recent reviews that discuss the significance of such near-critical
models are \cite{BG11} and \cite{MM12}. 
%The near-critical field can be shown to be massive --- i.e., have
%exponential decay of correlation --- by using FK percolation methods,
%but we will not do so in this paper.
%\margin{added sentence and two new refs.}
%We will use the following definition for a full-plane field $\Phi$
%to be {\it massive\/}: namely, that for any bounded test functions $f,g$
%of compact support such that $\Phi(f)$ and $\Phi(g)$ have finite second
%moments, 
%\[
%|{\text{Cov}} (f, T_y g)| \, \leq \, C \, e^{- \mu ||y||} 
%\]
%for some $C=C(f,g) \in (0,\infty)$ and $\mu > 0$ (not depending
%on $f,g$ and $y$). Here $T_y$ denotes translation of $f$ by $y \, \in \, \R^2$
%and $||y||$ denotes the Euclidean norm of $y$. We remark that such
%exponential decay for the two-point correlation of $\Phi$ 
%implies exponential decay of $n$-point correlations when $\Phi$
%satisfies the FKG inequalities --- see, e.g., \cite{Leb72} and
%\cite{Sim73}.

\begin{theorem}\label{th.NC}
Let us fix some constant $h>0$. 
Consider the Ising model on $a\, \Z^2$ at $\beta=\beta_c$ and with  vanishingly small external magnetic field equal to $a^{15/8}h$. Let $\Phi^{a,h}$ be the near-critical magnetization field in the plane defined, as in \cite{CGN1} (where $h = 0$), by 
\begin{align*}
\Phi^{a,h}:= \sum_{x\in a\Z^2} \, \delta_x\, \sigma_x\, a^{15/8}\,,
\end{align*}
where $\{\sigma_x\}_{x\in a\Z^2}$ is a realization of the above 
Ising model with external magnetic field equal to $h\, a^{15/8}$. 
Then, as the mesh $a\searrow 0$, the random distribution $\Phi^{a,h}$ 
converges in law to a 
%{\bf massive magnetization field} $\Phi^{\infty, h}$ under 
near-critical field $\Phi^{\infty, h}$ under
the topology of $\Sob^{-3}$ in the full plane defined in Section A.2 of \cite{CGN1}. 
\end{theorem}

The analogous statement in the case of a bounded smooth domain can be stated as follows.
%Let us point out that the corresponding statement in a bounded smooth domain $\Omega$ is much easier to obtain using the readily the results from \cite{CGN1}. We have
\begin{proposition}\label{pr.NC}
Let $\Omega$ be a bounded smooth domain of the plane with boundary conditions either $+,-$ or $\free$ and let $h>0$ be some positive constant.  
Then, with the obvious notation, $\Phi^{a,h}_\Omega$ converges in law to 
a field $\Phi^{\infty,h}_\Omega$ as $a\to 0$ under the topology of the Sobolev space $\Sob^{-3}(\Omega)$. 
\end{proposition}
This result is stated only as a proposition since as we shall see in 
Section \ref{s.NC}, it is follows almost readily from our previous 
work \cite{CGN1}. We will then prove Theorem \ref{th.NC} using 
this proposition by considering larger and larger domains $\Lambda_L$ 
and by showing that the near-critical fields do stabilize as $L\to \infty$. 
 The relation between $\Phi^{\infty, h}$ and $\Phi^{\infty,h}_\Omega$
to $\Phi^{\infty}$ and $\Phi^{\infty}_\Omega$ is discussed in Section~\ref{s.NC}.

%\begin{theorem}\label{th.NC}
%\bi
%\item[(i)] Let $\Omega$ be bounded smooth domain of the plane.
%\ei
%\end{theorem}

\subsection{Brief outline of proofs}

\bi

\item The proof of the tail behaviour given by Theorem \ref{th.tail} 
will be based on the study of the exponential moments of the magnetization $m$, i.e., 
on $\Eb{e^{t\, m}}$, with $t>0$ large. Theorem \ref{th.tail} will then follow from a 
specific {\bf Tauberian} theorem of Kasahara~\cite{Kasahara}.
One issue in this program is to show that the random variable $m$ indeed has 
exponential moments. This property was established in the first part of 
this series of papers, i.e. in \cite{CGN1} and the proof relied essentially 
on the  {\bf GHS}~inequality. The other difficulty is to adapt the classical arguments 
which lead to the existence of {\it free energies} to our present continuum setting, where one cannot use the standard trick of fixing the spins along dyadic squares in order to
 use subadditivity. To overcome this, one relies on RSW within thin long tubes. 

\item In our study of $|\EFK{}{\xi}{e^{i\, t\, m}}|$, we rely on the FK representation 
and we prove that with very high probability (of order $1-e^{-c |t|^{16/15}}$), 
one can find $O(1/\eps^2)$ mesoscopic squares of 
well-chosen size $\eps=\eps_t$ which contain an FK cluster of ``mass'' about $1/t$. 

\item For the proof of Theorem \ref{th.NC}, most of the non-trivial work 
is done in Lemma \ref{l.coupling} whose purpose is to prove that the law of the full-plane 
near-critical field $\Phi^{a,h}$ is very close in $\Sob^{-3}$ to the law of a large domain 
near-critical field $\Phi^{a,h}_{\Lambda_L}$. The technique used here is a coupling argument 
similar to the one used in Section 2 in \cite{CGN1} and which relies on the RSW Theorem from \cite{RSWfk}.
\ei

\medskip 
\ni
{\bf Acknowledgments.}
We wish to thank Hugo Duminil-Copin for his insights which lead to Remark~\ref{r.xi}.

\section{Tail behavior}\label{s.tail}

In this section, we shall prove Theorem \ref{th.tail}.

\subsection{Existence of exponential moments}

We will need the fact that the (continuum) magnetization $m$ has all exponential moments.
This property was proved in \cite{CGN1} and we provide below the corresponding statement.

\begin{proposition}[Proposition 3.5 and Corollary 3.8  in \cite{CGN1}]\label{pr.expmom}
For any boundary condition $\xi$ (either  $+$, $-$ or $\free$ boundary conditions) around $[0,1]^2$,
and for any $t\in \R$, if $m=m^\xi$ is the continuum magnetization of the unit square, then one has 
\bi
\item[(i)] $\Eb{e^{t\, m}} < \infty$.
\item[(ii)] Furthermore, as $a\to 0$, 
$\Eb{e^{t\, m^a}} \to \Eb{e^{t\, m}}$.
\ei

\end{proposition}

See \cite{CGN1} for the proof of this proposition which relies essentially on the {\bf GHS inequality} from \cite{GHS}.
%
%\begin{theorem}[GHS inequality, \cite{GHS}]\label{th.GHS}
%Let $G=(V,E)$ be a finite graph. Consider a ferromagnetic Ising model on this graph (i.e. the interactions between $i\sim j$ are positive)
%and assume furthermore that the external field $\mathbf{h} = (h_v)_{v\in V}$ (which may vary from one vertex to another) is positive. Under such general assumptions, one has for any vertices 
%$i,j,k\in V^3$:
%\[
%\langle \sigma_i\sigma_j \sigma_k\rangle - 
%\Bigl( 
%\langle \sigma_i\rangle\, \langle \sigma_j \sigma_k\rangle
%+\langle \sigma_j\rangle\, \langle \sigma_i \sigma_k\rangle
%+\langle \sigma_k\rangle\, \langle \sigma_i \sigma_j\rangle
%\Bigr)
%+2 \langle \sigma_i\rangle \langle \sigma_j\rangle \langle \sigma_k\rangle
%\leq 0 \,.
%\]
%\end{theorem}

\subsection{Asymptotic behavior of the moment generating function and scaling argument}

Since the exponential moments $\EFK{}{\xi}{e^{t\, m}}$ are well-defined, our next step is to study the behavior for large $t$ of the moment generating function
$t\mapsto \EFK{}{\xi}{e^{t\, m}}$. We will prove the following proposition.

\begin{proposition}\label{pr.asympLap}
There exists a {\bf universal constant} $b>0$ which does not depend on the boundary conditions $\xi$ around $[0,1]^2$ so that 
as $t\to \infty$:
\[
\log \EFK{}{\xi}{e^{t\, m}} \sim b \; t^{\frac {16}{15}}\,.
\]
\end{proposition}

Theorem \ref{th.tail} follows from the above proposition thanks to the following Tauberian Theorem by Kasahara. 

\begin{theorem}[Corollary 1 in \cite{Kasahara} ]\label{th.Taub}
For any random variable $X$ which has all its exponential moments, if there is an exponent $\alpha>1$ and a constant $b>0$ such that 
\begin{equation*}
\log \Eb{e^{t\, X}} \sim  b\, t^\alpha \,,
\end{equation*}
as $t\to \infty$, 
then the following holds for some explicit constant $c=c(b,\alpha)>0$:
\begin{equation*}
\log \Pb{X>x}\sim  - c\, x^{\frac 1 {1-1/\alpha}}\,,
\end{equation*}
as $x\to \infty$. 
\end{theorem}

\begin{remark}\label{}
In fact, this result is stated only for positive random variables in \cite{Kasahara} but it is very simple to extend it to any  real-valued random variable $X$. Let us sketch a short argument here.
Assume one has
\begin{equation}\label{}
\log \Eb{e^{t\, X}} \sim  b\, t^\alpha \,,
\end{equation}
as $t\to \infty$ for some $b>0$, then necessarily, $\Pb{X>0}$ has to be 
strictly positive. Now let $Y$ be the random variable  $X$ conditioned to be positive. 
It is easy to check that as $t\to \infty$,
$\log \Eb{e^{t\, Y}} \sim \log \Eb{e^{t\, X}}$.
One then concludes the argument by noticing that as $x\to \infty$, 
$\log \Pb{X>x} \sim \log \Pb{X>x \md X>0} = \log \Pb{Y>x}$. 
\end{remark}

\begin{remark}
Note that by a straightforward use of the exponential Chebyshev inequality, upper bounds on $\FK{}{\xi}{m>x}$ 
 can be directly recovered from Proposition \ref{pr.asympLap}.
\end{remark}

\medskip

\ni
{\bf Proof of Proposition \ref{pr.asympLap}:}

The main tools to prove the Proposition will be the scaling covariance property of the total magnetization $m$ which was proved in \cite{CGN1} (see Proposition \ref{pr.SC} below) as well as Theorem \ref{th.FE} below which in some sense defines a {\bf free energy} for our limiting magnetization field.
Let us first state these two results.

\begin{proposition}[Scaling covariance of $m$, Corollary 5.2 in \cite{CGN1}]\label{pr.SC}
Let $m=m^\xi$ be the scaling limit of the renormalized magnetization in the square (i.e., $m = \<{\Phi^\infty, 1_{[0,1]^2}}$), with boundary conditions $\xi$ being either $+,-$ or $\free$.
For any $\lambda>0$, let $m_\lambda=m_\lambda^\xi$ be the scaling limit of the renormalized magnetization in the square $[0,\lambda]^2$ with the same boundary conditions $\xi$. Then one has the following identity in law:
\begin{equation}\label{}
m_\lambda \overset{(d)}{=}  \lambda^{15/8}\, m\,.
\end{equation}
\end{proposition}

\begin{theorem}[Existence of free energy]\label{th.FE}
For any $L>0$ and any boundary conditions 
$\xi$  (made of finitely many  $+$, $-$ or $\free$ arcs) around $[0,L]^2$,
let $f_L^{\xi}(t):= \frac 1 {L^2}  \log \EFK{}{\xi}{e^{t\, m_L}}$. 

There is a universal constant $b>0$, which does not depend on the boundary conditions $\xi$, such that 
for any $t\in \R$
\[
f_L^{\xi}(t):= \frac 1 {L^2}  \log \EFK{}{\xi}{e^{t\, m_L}}  \underset{L\to \infty}{\longrightarrow} b\, |t|^{16/15}\,.
\]
\end{theorem}

With these two ingredients, it is easy to conclude the proof of Proposition \ref{pr.asympLap}. Indeed if $\lambda_t:=t^{8/15}$,
then one has:
\begin{align}\label{}
\log \EFK{}{\pm}{e^{t\, m}} & =  \log \EFK{}{\pm}{e^{ m_{\lambda_t}}}  \;\;\text{ using Proposition \ref{pr.SC}}\\
&= t^{16/15}  \Bigl( \frac{1} {\lambda_t^2}  \log \EFK{}{\pm}{e^{ m_{\lambda_t}}}  \Bigr) \\
&\sim t^{16/15} \, b\,,
\end{align}
as $t \to \infty$.
Other boundary conditions are handled by noting that $\EFK{}{\xi}{e^{t\, m}}$ is squeezed
between the $+$ and $-$ cases by the FKG inequalities. 
\QED

\begin{remark}\label{}
Note that we did not need the full strength of Theorem \ref{th.FE}, only the case $t=1$. 
Nevertheless, since Theorem \ref{th.FE} is interesting in it own right, we prove it for all $t\in \R$ (which will result in a slight repetition of the above scaling argument in the proof of Lemma \ref{l.final}). 
\end{remark}

\begin{remark}\label{}
It is tempting to compare the above {\bf free energy} with the classical one coming from the discrete system, 
i.e., defined as 
\begin{equation}\label{}
F(t):= \lim_{N\to \infty} \frac 1 {N^2} \EFK{}{+}{e^{t\, \sum_{x\in \Lambda_N} \sigma_x}}\,,
\end{equation}
but it is easy to see that they must be different, since clearly $F(t) \le |t|$ for any $t\in \R$. 
On the other hand, they behave essentially the same for small $t$
as follows from the results of~\cite{CGN2}.
\end{remark}

\subsection{Free energy estimates}

The purpose of this section is to prove Theorem \ref{th.FE} on the free energy of $\Phi^\infty$.
The proof of this theorem will be divided into several steps as follows. 
First, we will show in Lemma \ref{l.dyadic} that for any $t \geq 0$, $f^+_L(t)$ and $f^-_L(t)$ have limits along dyadic scales $L_k=2^k$, respectively denoted by $f^+(t)$ and $f^-(t)$. 
Then, in Lemma \ref{l.trueL}, we will show that 
\begin{equation*}
\begin{cases}
\limsup_{L\to \infty}f^+_L(t) &= f^+(t) \\
\liminf_{L\to \infty} f^-_L(t) &=f^-(t)\,.
\end{cases}
\end{equation*} 
In Lemma \ref{l.+=-}, we will prove that $f^+(t)=f^-(t)=f(t)$ for any $t\geq 0$. 
Finally Lemma \ref{l.final} will identify the limit $f(t)$ to be exactly $b\, |t|^{16/15}$ for all $t\in \R$, thus concluding the proof of Theorem \ref{th.FE}.
The main difficulty in this last lemma will be to show that the constant $b$ is positive.

We will first list these lemmas and then proceed with their proofs. 
Let us point out that some of the proofs below follow the standard arguments to prove that a {\it free energy} is well defined.
Nevertheless, they turn out to be slightly more involved  here since we are working with the continuum limit and therefore all the classical arguments based, for example, on counting the number of {\it lattice sites} on the boundary are no longer valid here. (Only the proof of Lemma \ref{l.dyadic} follows exactly the classical scheme).

\begin{lemma}\label{l.dyadic}
For any $t \geq 0$, and any $k\geq 1$,
\begin{equation*}
\begin{cases}
f^+_{2^{k+1}}(t) \le f^+_{2^k}(t) \, , \\
f^-_{2^{k+1}}(t) \ge f^-_{2^k}(t)\,.
\end{cases}
\end{equation*}
In particular, the sequences $f^\pm_{2^k}(t)$ converge as $k \to \infty$ and we will denote respectively their limits by $f^\pm(t)$. 
\end{lemma}

\begin{lemma}\label{l.trueL}
For any $t\geq 0$, we have 
\begin{equation*}
\begin{cases}
\limsup_{L\to \infty}f^+_L(t) &= f^+(t) \, , \\
\liminf_{L\to \infty} f^-_L(t) &=f^-(t) \, . 
\end{cases}
\end{equation*} 
\end{lemma}

\begin{lemma}\label{l.+=-}
For any $t\geq 0$, we have 
\begin{equation*}
f^+(t)=f^-(t) = f(t)\,.
\end{equation*}
\end{lemma}

\begin{lemma}\label{l.final}
There exists a universal constant $b>0$ such that for any boundary conditions $\xi$, we have 
\begin{equation*}
f(t):= \lim_{L\to \infty} \frac 1 {L^2} \EFK{}{\xi} {e^{t\, m_L}}  = b \, |t|^{16/15}\,,
\end{equation*}
for all $t\in \R$. 
\end{lemma}

\ni
{\bf Proof of Lemma \ref{l.dyadic}:}
Let us consider the case of $+$ boundary conditions;
the $-$~case is similar.
From Proposition \ref{pr.expmom} (ii), we know that for any $L>0$, 
\begin{equation*}
\EFK{}{+}{e^{t\,m_L}} = \lim_{a\to0} \EFK{}{+}{e^{t\, m^a_L}}\,.
\end{equation*}
Now, for any $k\in \N^+$, it is easy to check (by breaking the domain $[0,2^{k+1}]^2$ into $4$
squares with $+$ boundary conditions and using FKG) that for suitable choices of the mesh size $a$
(i.e. $a$ such that $2^k\in a\Z$), then  
\begin{align*}\label{}
\EFK{}{+}{e^{t\, m^a_{2^{k+1}}}} & \le \EFK{}{+}{e^{t\, m^a_{2^k}}}^{4}\,.
\end{align*}
Taking the limit $a\to 0$, we get that 
\begin{align*}\label{}
\EFK{}{+}{e^{t\, m_{2^{k+1}}}} & \le \EFK{}{+}{e^{t\, m_{2^k}}}^{4}\,,
\end{align*}
which implies $f_{2^{k+1}}^+(t) \le f_{2^k}^+(t)$. As pointed out above, this proof matches exactly the standard proof in the discrete setup.
\QED

\ni
{\bf Proof of Lemma \ref{l.trueL}:}

\ni

We only consider the case of $+$ boundary conditions and we will fix some $t\geq 0$ (the case of minus boundary conditions is handled in the same fashion). 
Let us also fix some integer $k_0 \geq 1$. We wish to show that   $\limsup_{L\to \infty} f^+_L(t) \le f^+_{2^{k_0}}(t)$.

%Let us fix some small $\eps>0$. Let $k_0\geq 1$ be large enough so that for all $k\geq k_0$, 
%\[
%|f^+_{2^k}(t)-f^+(t)|<\eps\,.
%\]

\begin{figure}[!htp]
\begin{center}
\includegraphics[width=0.5\textwidth]{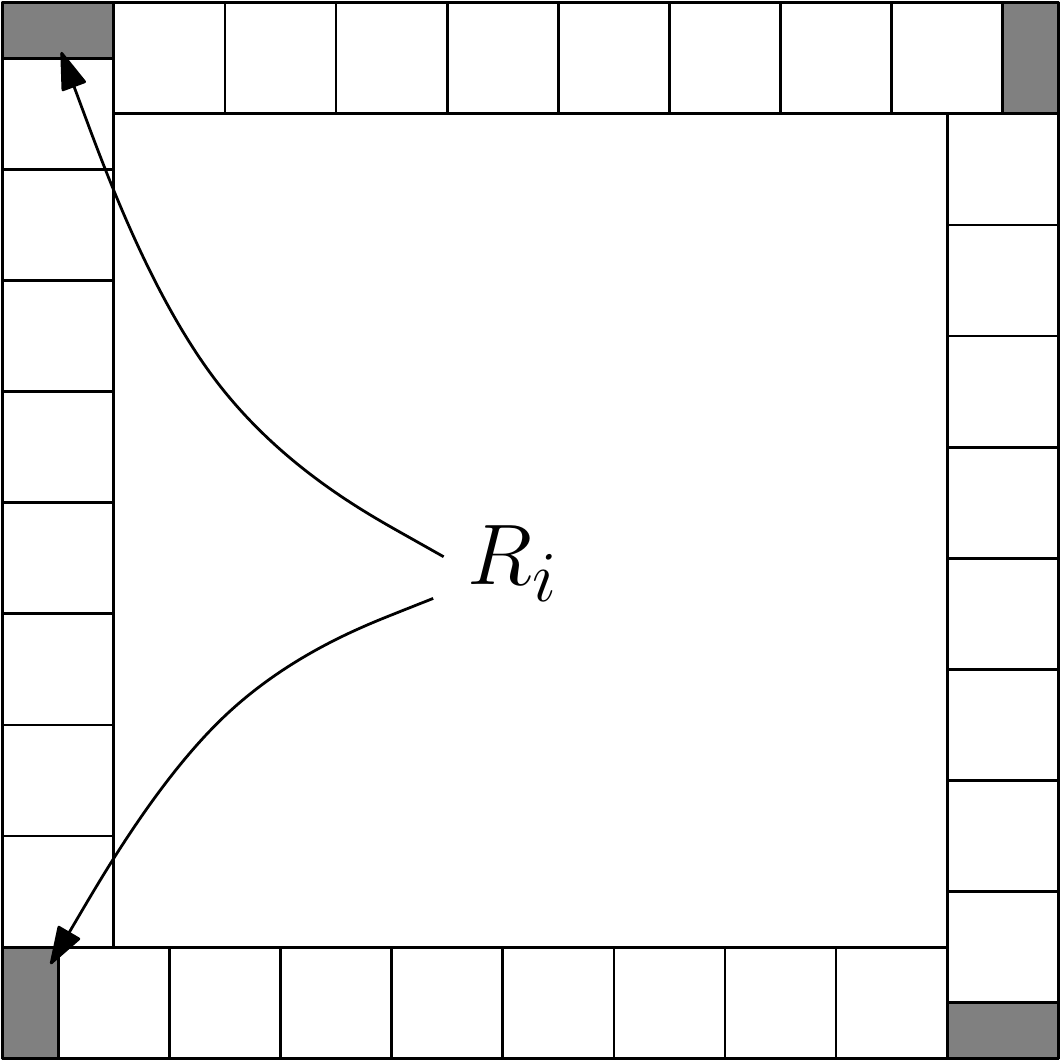}
\end{center}
\caption{The annulus $A$ and the rectangles $R_1,\ldots,R_4$.}\label{f.rectangle}
\end{figure}

For $L>0$ large enough, let $M\geq 1$ be such that $L= M\, 2^{k_0} + 2 K$, with $K\in [2^{k_0-1},2^{k_0})$. 
Divide the domain $[0,L]^2$ into the inside square $Q:=[K,L-K]^2$ and the annulus $A:=[0,L]^2\setminus Q$. 
Then, as in the proof of the above Lemma, we have 
\begin{align}\label{e.controlK}
\EFK{}{+}{e^{t\, m_{L}}} & \le \EFK{}{+}{e^{t\, m_{2^{k_0}}}}^{M^2} \EFK{}{+}{e^{t\, m_A}}\,,
\end{align}
where $m_A$ denotes the magnetization in the annulus $A$ with $+$ boundary conditions on its inner and outer boundaries.
One can split this annulus into %$O(L)$ 
 a number ($\leq 4L/K$) of
squares of side-length $K$ plus possibly 4 identical rectangles (up to a rotation) with one side of length $K$ and the other side of shorter length $\tilde K$---see Figure \ref{f.rectangle}. Call $R_1, \ldots, R_4$ those rectangles and let $\mathcal{R}$ be the family of possible shapes they can have.
Then, we have
\begin{align}\label{e.split}
\EFK{}{+}{e^{t\, m_A}} & \le \EFK{}{+}{e^{t\,m_K}}^{4L/K} \sup_{R\in \mathcal{R}} \EFK{}{+}{e^{t\, m_R}}^4 \nn \\
& \le \EFK{}{+}{e^{t\,m_K}}^{8(M+1)}  \sup_{R\in \mathcal{R}} \EFK{}{+}{e^{t\, m_R}}^4\,.
\end{align}

Now for any rectangle $R=[0,\tilde K] \times [0,K]$ with $0<\tilde K \le K$, one has  
\begin{align*}\label{}
 \EFK{}{+}{e^{t\, m_R}} & \le  \lim_{a\to0} \EFK{}{+}{e^{t\, m^a_R}} \\
& \le \lim_{a\to 0}  \exp \Bigl( t \,\EFK{}{+}{m^a_R} + \frac {t^2} 2 \EFK{}{+}{(m^a_R - \langle m^a_R\rangle)^2}\Bigr) \,,
\end{align*}
using the GHS inequality (see Theorem 3.6 and Corollary 3.7 of \cite{CGN1}).
As in the Appendix B of \cite{CGN1}, it is easy to check that 
\begin{equation*}
\sup_{R\in \mathcal{R}} \Bigl( \limsup_{a\to 0} \left\{ \EFK{}{+}{m^a_R} + \EFK{}{+}{(m^a_R)^2}\right\} \Bigr)< \infty\,,
\end{equation*}
which thus  implies 
\begin{equation}\label{e.uniform}
\sup_{R\in \mathcal{R}} \EFK{}{+}{e^{t\, m_R}} < \infty\,.
\end{equation}

In the same fashion, we have that 
\begin{equation}\label{e.uniform2}
\sup_{K \in [2^{k_0-1}, 2^{k_0}) } \EFK{}{+}{e^{t\, m_K}}< \infty\,.
\end{equation}

Plugging  the previous estimates into ~\eqref{e.controlK}, we obtain 
\begin{align*}
\frac 1 {L^2} \log \EFK{}{+}{e^{t\, m_L}} &  \le  \frac {M^2} {L^2} \log \EFK{}{+}{e^{t\, m_{2^{k_0}}}} + \frac{ 8M+8}{L^2}  \sup_{K} \log \EFK{}{+}{e^{t\,m_K}} + \frac 4 {L^2}  \sup_{R\in \mathcal{R}} \log \EFK{}{+}{e^{t\, m_R}}\,.
\end{align*}

By letting $L,M \to \infty$, the last two terms tend
to zero, while the first one converges to $f^+_{2^{k_0}}(t)$, which ends the proof of the lemma.
\QED.

\ni
{\bf Proof of Lemma \ref{l.+=-}: }

\ni
It is clear, by monotonicity, that for any $t\geq 0$, $f^-(t) \le f^+(t)$. 
Let us then show the reverse inequality. 
We will in fact compare the plus boundary conditions with free boundary conditions showing with the obvious notation
that $f^\free(t) \geq f^+(t)$. Since the same proof allows us to show that $f^\free(t) \le f^-(t)$,
this is enough to conclude the proof.

We wish to show that 
\begin{align*}\label{}
f^\free(t):= \liminf_{k\to \infty} 2^{-2k} \log \EFK{}{\free}{e^{t\, m_{2^k}}}  \geq f^+(t)=  \lim_{k\to \infty} 2^{-2k} \log \EFK{}{+}{e^{t\, m_{2^k}}}\,.
\end{align*}
Note that we used $\liminf$ to define $f^\free$ here since we have not proved (yet) that the limit exists in the case of $\free$ boundary conditions and $\liminf$ is the worst possible case here.

Let us fix some small dyadic $\eps=2^{-k_0}>0$.
For any $L>10$, let $R=R_{L,\eps}$ be the event that there is a + cluster in the annulus $a\Z^2 \cap [0,L]^2 \setminus [\eps L , (1-\eps)L]^2$.  From the RSW Theorem in \cite{RSWfk}, we have that 
\begin{equation*}
H:=\inf_{L>10, a<1}\FK{}{\free}{R_{L,\eps}}>0\,. 
\end{equation*}
Recall furthermore that for any $L>10$: 
\begin{equation*}
\frac 1 {L^2} \log \EFK{}{\free}{e^{t\, m_L}} = \lim_{a\to 0} \frac 1 {L^2} \log \EFK{}{\free}{e^{t\, m_L^a}}\,.
\end{equation*}

We have that 

\begin{align*}\label{}
\liminf_{L\to \infty} \frac 1 {L^2} \log \EFK{}{\free}{e^{t\, m_L}} &
 \geq \liminf_{L\to \infty} \lim_{a\to 0} \frac 1 {L^2} \log \EFK{}{\free}{ 1_{R}\, e^{t\, m_L^a}} \\
&\geq  \liminf_{L \to \infty}  \lim_{a\to 0} \frac 1 {L^2} \log \EFK{}{\free}{ e^{t\, m_L^a} \, \md R}   + \lim_{L\to \infty} \frac 1 {L^2} \log  H \\
&=   \liminf_{L \to \infty} \lim_{a\to 0}  \frac 1 {L^2} \log \EFK{}{\free}{ e^{t\, m_L^a} \, \md R}
\end{align*}

For each dyadic $L=L_k=2^k>10$, let us divide the square $[0,L]^2$ into the annulus $A=A_L = [0,L]^2 \setminus [\eps L , (1-\eps)L]^2$ 
and the inside square $Q=Q_L= [\eps L , (1-\eps)L]^2$. As such and with the obvious notation, we will decompose the magnetization $m_L^a$ into 
\begin{equation}\label{}
m_L^a = m_A^a + m_Q^a\,.
\end{equation}
Furthermore, we will denote by $\F_Q$ the filtration generated by the spins in $a\Z^2\cap Q$. By conditioning furthermore on $\F_Q$, we get 
\begin{align}\label{e.FQ}
\liminf_{k\to \infty} \frac 1 {L_k^2} \log \EFK{}{\free}{e^{t\, m_{L_k}}} & \geq 
\liminf_{k \to \infty}  \lim_{a\to 0} 
\frac 1 {L_k^2} \log \EFK{}{\free}{ e^{t\, m_Q^a} \, \EFK{}{\free}{e^{t\, m_A^a} \md \F_Q, R} \, \md R} \, .
\end{align}

Let us first show the following lemma.
\begin{lemma}\label{l.FQ}
 There is a function $\alpha=\alpha(\eps)>0$ such that uniformly in $L=2^k>10$ %, in $a<1$
and on the configuration of spins $\sigma_Q$ inside $Q$, one has 
\begin{equation}\label{e.FQ2}
\lim_{a\to 0}\EFK{}{\free}{e^{t\, m_A^a} \md \F_Q, R} \geq \alpha(\eps)\,.
\end{equation}
\end{lemma}

\ni
{\bf Proof:}

To prove the lemma, notice that by our choice of $\eps$, the annulus $A$ can be divided into $4 (2^{k_0} -1)$ exact squares of side-length $2^{k-k_0}$ (as in Figure \ref{f.rectangle} except there are no thin rectangles there) and we have the bound 
\begin{align}\label{}
\lim_{a\to 0} \EFK{}{\free}{e^{t\, m_A^a} \md {\cal F}_Q,R} & \geq \lim_{a\to 0 }\EFK{}{-}{e^{t\, m_{2^{k-k_0}}^a}}^{4 (2^{k_0} -1)}
\;\;\text{ using FKG} \\
&\geq \lim_{a\to 0} \FK{}{-}{m_{2^{k-k_0}}^a > 0}^{4 (2^{k_0} -1)} \\
&\geq \FK{}{-}{m_{2^{k-k_0}} > 0}^{4 (2^{k_0} -1)}  \\
& = \FK{}{-}{m_{[0,1]^2} > 0}^{4 (2^{k_0} -1)}  \geq  \FK{}{-}{m_{[0,1]^2} > 0}^{4/\eps}\,,
\end{align}
%\margin{please check the use of FKG is fine here}
where in the last line, we relied on the scaling covariance property given by Proposition \ref{pr.SC}.
\QED

We conclude the proof of Lemma \ref{l.FQ} by relying on the following easy lemma.
\begin{lemma}\label{l.easy}
There is a constant $c>0$ such that 
\begin{equation*}
\FK{}{-}{m_{[0,1]^2} > 0}>c\,.
\end{equation*}
\end{lemma}

%\ni
%{\bf Sketch of proof of the Lemma:}
This lemma can be proved for example by using the FK representation of the field 
$\Phi^\infty$ from \cite{CGN1} and the fact that small FK clusters contribute little 
to the total magnetization $m_{[0,1]^2}$ (see for example 
%Theorem 3.9 of \cite{CGN1} \margin{Shouldn't we refer to Eqs. (2.8)-(2.11) 
Equations (2.8)--(2.11) of \cite{CGN1}).  
%of \cite{CGN1}?}). 
Hence this ends the proof of  Lemma \ref{l.FQ} with $\alpha(\eps):=c^{4/\eps}$. \QED

Plugging~\eqref{e.FQ2} into ~\eqref{e.FQ} gives us 
\begin{align}\label{e.FQ3}
\liminf_{k\to \infty} \frac 1 {L_k^2} \log \EFK{}{\free}{e^{t\, m_{L_k}}} & \geq 
\liminf_{k \to \infty}   \lim_{a\to 0} \frac 1 {L_k^2} \log \EFK{}{\free}{ e^{t\, m_Q^a} \md R} \, .
\end{align}

Now, by FKG it is clear that 
\begin{align} \label{new-label}
\lim_{a\to 0}  \EFK{}{\free}{ e^{t\, m_Q^a} \md R} & \geq \lim_{a \to 0}\EFK{}{+}{e^{t\, m_Q^a}} \\
& = \EFK{}{+}{e^{t\, m_Q}}\,, \nonumber
\end{align}
where in the latter expectations, the $+$ boundary conditions are around $[0,L_k]^2$ and hence are further away from the domain $Q=Q_{L_k}$. 

To conclude the proof of Lemma \ref{l.+=-} we still need to compare $\EFK{}{+}{e^{t\, m_Q}}$ with $\EFK{}{+}{e^{t\, m_L}}$. This is done by the following lemma.
\begin{lemma}\label{l.end}
There is a function $\eta(x)>0$ satisfying $\eta(x)\to 0$
as $x\to 0$, and such that for any $\eps=2^{-k_0}$, one has, with the above notation,

\begin{equation}\label{}
\lim_{k\to \infty} \frac 1 {L_k^2} \log \EFK{}{+}{e^{t\, m_Q}} \geq  \lim_{k\to \infty} \frac 1 {L_k^2} \log \EFK{}{+}{e^{t\, m_{L_k}}} - \eta(\eps)\,.
\end{equation}
\end{lemma}

\ni
{\bf Proof:}%of Lemma \ref{l.end}:}

As in the proof of Lemma \ref{l.trueL}, and dividing $[0,L_k]^2$ as above,
%using the same division as the one above,
we have 
\begin{align}\label{}
\EFK{}{+}{e^{t\, m_{L_k}}} & =  \EFK{}{+}{e^{t\, m_{Q} + t\, m_A }} \\ 
& \le \EFK{}{+}{e^{t\, m_Q}} \EFK{}{+}{e^{t\, m_{2^{k-k_0}}}}^{4(2^{k_0}-1)}\,,
\end{align}
where, as above, the boundary conditions in the expectation $\EFK{}{+}{e^{t\, m_Q}}$ are meant to be 
around the larger square $[0,L_k]^2$. 
Now, we have 
\begin{align}\label{}
\frac{1} {L_k^2} \log \EFK{}{+}{e^{t\, m_{2^{k-k_0}}}}^{4(2^{k_0}-1)} &
= \frac {4(2^{k_0}-1)}{2^{2 k_0}} \frac{1}{2^{2(k-k_0)}} \log \EFK{}{+}{e^{t\, m_{2^{k-k_0}}}} \\
&\le 4\eps \; f^+_{L_{k-k_0}}(t)\,.
\end{align}
Letting $k\to \infty$, we obtain
\begin{equation}\label{last-eq}
\lim_{k\to \infty} \frac 1 {L_k^2} \log \EFK{}{+}{e^{t\, m_Q}} \geq f^+(t) - 4\eps f^+(t)\,.
\end{equation}
This ends the proof of Lemma \ref{l.end}. \QED

To conclude the proof of Lemma \ref{l.+=-}, we plug \eqref{last-eq} into ~\eqref{e.FQ3} and obtain,
using \eqref{new-label}, 
\begin{equation*}\label{}
\liminf_{k\to \infty} \frac 1 {L_k^2} \log \EFK{}{\free}{e^{t\, m_{L_k}}} \geq f^+(t) - 4\eps f^+(t)\,,
\end{equation*}
for any value of $\eps=2^{-k_0}>0$. Hence, we have that 
\begin{equation}\label{}
\liminf_{k\to \infty} \frac 1 {L_k^2} \log \EFK{}{\free}{e^{t\, m_{L_k}}} \geq f^+(t)\,,
\end{equation}
which thus implies
\begin{equation}\label{}
 f^\free(t) = \liminf_{k\to \infty} \frac 1 {L_k^2} \log \EFK{}{\free}{e^{t\, m_{L_k}}}   = f^+(t)\,.
\end{equation}
 \QED

\ni
{\bf Proof of Lemma \ref{l.final}:}

As in the proof of Proposition \ref{pr.asympLap}, using the scaling covariance given by Proposition \ref{pr.SC}, we have that for any $L>0$ and any $t>0$ and for, say, $+$ boundary conditions,
\begin{equation*}
t\, m_L  \overset{(d)}{=} m_{L\, t^{8/15}}\,.
\end{equation*}
This implies 
\begin{align*}
f(t)=f^+(t) &= \lim_{L\to \infty} \frac 1 {L^2} \log \EFK{}{+}{e^{t\, m_L}} \\
&= \lim_{L \to \infty} \frac 1 {L^2}  \log \EFK{}{+}{e^{ m_{L\, t^{8/15}}}} \\
&= t^{16/15}\,  \lim_{\bar L \to \infty} \frac 1 {\bar L^2}  \log \EFK{}{+}{e^{ m_{\bar L}}} \\
&= f(1) \, t^{16/15}\,.
\end{align*}
To conclude the proof of the lemma when $t>0$, it remains to show that the quantity (with $L_k=2^k$)
\begin{equation*}
f(1)=f^+(1)=\lim_{k\to \infty} \frac 1 {L_k^2} \log \EFK{}{+}{e^{m_{L_k}}}
\end{equation*}
is strictly positive. %($>0$).

To see this, let us first denote by $M_L$ the magnetization
$\Phi^\infty(1_{[0,L]^2})$ in $[0,L]^2$ of the full-plane field
$\Phi^{\infty}$. By the results of~\cite{CGN1}, for any
$L \in (0,\infty)$, $M_L$ has zero mean and variance in $(0,\infty)$.
Then by a few uses of the FKG inequalities, we have
\begin{align*}
%\EFK{}{+}{e^{m_{2^k}}} &\geq \EFK{}{-}{e^{m_1}}^{L_k^2} \\
%&\geq \FK{}{-}{m>0}^{L_k^2}\,,
\EFK{}{+}{e^{m_{2^k}}} &\geq \EFK{}{}{e^{M_{2^k}}} \\
&\geq (\EFK{}{}{e^{M_1}})^{L_k^2} \\
&\geq (1 + \EFK{}{}{M_1^{\,2}})^{L_k^2} \, ,
\end{align*}
so that $f^+(1) \geq \log (1 + \EFK{}{}{M_1^{\,2}})\,> \,0$.
\QED

%This can be seen, for example using the lower bound  
%\begin{align*}
%\EFK{}{+}{e^{m_{2^k}}} &\geq \EFK{}{-}{e^{m_1}}^{L_k^2} \\
%&\geq \FK{}{-}{m>0}^{L_k^2}\,,
%\end{align*}
%and relying on Lemma \ref{l.easy}. 
%Finally the case $t=0$ in Lemma \ref{l.final} is clear and the case $t<0$ is done by noticing the fact that 
%when $t<0$, we have 
%\begin{equation*}
%\begin{cases}
%\EFK{}{+}{e^{t\, m}} = \EFK{}{-}{e^{|t|\, m}} \\
%\EFK{}{-}{e^{t\, m}} = \EFK{}{+}{e^{|t|\, m}} \, .
%\end{cases}
%\end{equation*}
%\QED

%\note{Reading note 21/06/13}
%{This whole section is done and seems correct (FKG checked and so on)}

\section{Analyticity of the probability density function of $\magn$}\label{s.analyt}

%Let us fix the boundary condition $\xi$ around $[0,1]^2$ to be either $+, -$ or $\free$.
In this section, we shall prove Theorem \ref{th.analyt}.
First of all, by the convergence in law of $m^a$ towards $m$, we have as $a\to 0$:
%$\Eb{e^{t\, m^a}} \to \Eb{e^{t\, m}}$ when $a\to 0$ (see item $(ii)$ in Proposition \ref{pr.expmom}), it is easy to check that the following holds as $a\to 0$:
\begin{align}\label{}
\EFK{}{\xi}{e^{i\, t\, m^a}} \to \EFK{}{\xi}{e^{i\, t\, m}}\,.
\end{align}
It is thus sufficient to prove that there exists a constant $c>0$ which is such that, for any $t\in \R$,
\[
\limsup_{a\to 0} |\EFK{}{\xi}{e^{i\, t\, m^a}}| \le e^{- c \, |t|^{\frac {16}{15}}}\,.
\]

To prove this, we will rely on the FK representation of the Ising model in $a\Z^2 \cap [0,1]^2$ endowed with its boundary conditions $\xi\in \{+,-,\free\}$. Let us assume that $\xi=+$. (The case of free boundary conditions is even easier). We can write  
\begin{align}\label{e.prod}
|\EFK{}{+}{e^{i\, t \, m^a}}| & = |\EFK{}{\fk}{  e^{i t \mathcal{A^+}}\, \prod_{C_i} \frac 1 2 (e^{i t \mathcal{A}_i} + e^{-i t \mathcal{A}_i}) } | \nn \\
&=  |\EFK{}{\fk}{  e^{i t \mathcal{A^+}}\, \prod_{C_i}  \cos \, t \mathcal{A}_i  }| \,,
\end{align}
where $\{C_i\}_i$ denotes the collection of clusters that do not intersect the boundary and $C^+$ is the cluster that intersects the boundary. Furthermore, we let $\mathcal{A}_i=\mathcal{A}_i^a=\mathcal{A}_i^a(C_i)$ stand for the renormalized areas of the cluster $C_i$, and $\mathcal{A}^+$ for the renormalized area of the cluster $C^+$. Our strategy, in order to obtain an upper bound for \eqref{e.prod}, is to show that with high probability, there are many
clusters in $\{ C_i \}_i$ with a renormalized area of order $1/t$. 

We will rely on the following lemma:
\begin{lemma}\label{l.area}
There exist constants $c\in(0,1)$ and $M>1$ such that 
for any $0<\eps <1/10$ and any $\eps$-square $Q$ inside $[0,1]^2$, %at distance at least $\eps$ from $\p [0,1]^2$,
uniformly as $a\to 0$, and uniformly on the FK configuration outside of $Q$, with (conditional) probability at least $c>0$, one can find at least one FK cluster $C$ inside $Q$ that does not intersect $\p Q$ and such that its renormalized area lies in the interval $[\eps^{15/8}/M, M\eps^{15/8}]$.  
\end{lemma}

\ni
{\bf Proof:}

Let $Q$ be an $\eps$-square inside $[0,1]^2$ and let $\omega_a$ be any FK-configuration outside $Q$. 
Let $A_1$ be the annulus $Q \setminus 7/8 Q$, $A_2$ the annulus $7/8 Q \setminus 3/4 Q$ and $A_3$ the annulus $3/4 Q \setminus 1/2 Q$.
Let us introduce the following events: let $D_1$ be the event that there is a {\bf dual} circuit in the annulus $A_1$ and let $O_2$ and $O_3$ be the events that there is an {\bf open} circuit around each annuli $A_2$ and $A_3$. 
 Using the RSW Theorem from \cite{RSWfk} for a free boundary condition, one has that there is a constant $b>0$ such that uniformly on the outside configuration $\omega_a$, one has 
 \begin{align}\label{}
\FK{}{+}{D_1, O_2,O_3 \md \omega_a} > b\,.
\end{align}
Now let $X=X_a$ be the number of points in $a\Z^2 \cap 1/2 Q$ which are connected via an open-arm to $\p (\frac 3 4 Q)$. Then using similar computations as in Proposition B.2 in \cite{CGN1} or in Lemma 3.1. in \cite{CGN2}, one can find a constant $C>0$ such that 
\begin{align}\label{}
\begin{cases}
\EFK{}{+}{X_a \md D_1, O_2, O_3,\omega_a}\geq (\eps/a)^{15/8}/C \\
\EFK{}{+}{X_a^2\md D_1, O_2, O_3,\omega_a}\leq C\, (\eps/a)^{15/4} \, .
\end{cases}
\end{align}

By a standard second-moment argument, and using the fact that all points 
counted in $X_a$ belong to the same cluster (thanks to $O_3$), 
one obtains that with positive conditional probability, one can find a cluster $C$ 
which does not intersect $\p Q$ and whose renormalized mass is larger than 
$1/M \eps^{15/8}$. (Note that the event $O_2$  is there to ensure 
some positive information inside $3/4 Q$.)

It remains to prove an upper bound. In the same way as $X_a$ is smaller than the actual number of points in the open cluster we are interested in, one can also introduce $\tilde X_a$ to be the number of points inside the whole square $Q$ which are connected to the boundary $\p (3/4 Q)$. This random variable dominates the size of the cluster we are interested in. It is enough to control its expectation and it is easy to see that, for a well-chosen constant $C>0$, one has 
\begin{align*}\label{}
\EFK{}{+}{\tilde X_a\md \omega_a} & \leq \EFK{}{+}{\tilde X_a \md \wired \, \p Q}  \\
& \leq C\, (\eps/a)^{15/8}\,.
\end{align*}
Since $\FK{}{+}{D_1, O_2,O_3 \md \omega_a} > b$, this implies \begin{align*}\label{}
\FK{}{+}{\tilde X_a \geq M (\eps/a)^{15/8} \md D_1, O_2,O_3, \omega_a} \leq \frac 1 M  \frac C b \,.
\end{align*}
By choosing $M$ large enough (so that the conditional probabilities of lower bound and upper bound don't add up to something larger than one),  one concludes the proof of the lemma. \qed

\ni
{\bf Proof of  Theorem \ref{th.analyt}:}
\ni

For any $|t|>100$, choose $\eps=\eps_t$ so that $M \eps^{15/8}=\frac 1 { |t|}$ (we use the constants from Lemma \ref{l.area}). Use a tiling of the square $[0,1]^2$ so that one has $\frac{1}{2}\eps^{-2}$ disjoint $\eps$-squares $Q$.
%as in Lemma \ref{l.area}.
Recall from that lemma that for each such square $Q$, the probability that one has 
a cluster inside $Q$ with renormalized area in $[(1/M) \eps^{15/8}, M \eps^{15/8}]$ 
is larger than $c>0$ uniformly on what may happen outside of $Q$. We thus expect that at least about $\frac c 2 \eps^{-2}$ squares $Q$ will contain such a cluster. Let $G$ be the event that 
at last $\frac c 4 \eps^{-2}$ squares $Q$ have a cluster with renormalized area in $[(1/M) \eps^{15/8}, M \eps^{15/8}]$. Then, by a classical Hoeffding inequality one has that 

\begin{align}\label{e.Gc}
\FK{}{+}{G^c} \leq e^{- d\, \eps_t^{-2}} = e^{- d M^{16/15}\, |t|^{16/15}}\,,
\end{align}
for some universal constant $d>0$. Now, on the event $G$, we have 
\begin{align*}\label{}
 |e^{i t \mathcal{A^+}}\, \prod_{C_i}  \cos \, t \mathcal{A}_i  | & \leq [\cos 1/M^2 ]^{(c/4)\,  \eps_t^{-2}} \\ 
 & = [\cos 1/M^2 ]^{(c/4) \, M^{16/15}\, |t|^{16/15}} \\ 
 & \leq e^{- \tilde{c} \, |t|^{16/15}}\,,
\end{align*}
for some well-chosen constant $\tilde{c}>0$. Combining the above estimate with equations~\eqref{e.prod} and~\eqref{e.Gc}, we thus end the proof of Theorem \ref{th.analyt} with a possibly smaller value of $\tilde{c}>0$ (due to $\FK{}{+}{G^c}$ as well as to the region $|t|\in[0,100]$). \qed

\begin{remark}\label{r.xi}

As suggested after Theorem \ref{th.analyt}, it should be possible to extend the above proof to basically any boundary conditions $\xi$ (with the only constraint that one can prove a scaling limit result for $m^a$ as in \cite{CGN1}). For example, if $\xi$ is made of a finite combination of $+,-,\free$ arcs, this is handled in \cite{CGN1}. In this latter case, one would rely on the following extension of~\eqref{e.prod}: 
  
\begin{align*}
|\EFK{}{\xi}{e^{i\, t \, m^a}}| & = |\EFK{}{\fk, \xi}{ \prod_{C_i} \frac 1 2 (e^{i t \mathcal{A}_i} - e^{-i t \mathcal{A}_i}) \, \prod_{C_k^+} e^{i t \mathcal{A}_k}\, \prod_{C_l^-} e^{- i t \mathcal{A}_l} }| \\
& \leq  \EFK{}{\fk,\xi}{ \prod_{C_i}  | \cos \, t \mathcal{A}_i  | } \,.
\end{align*}

The additional difficulty when $\xi$ is a general boundary condition lies in the FK-representation of the associated Ising model. Indeed, general boundary conditions induce negative information in the bulk (since the FK configuration is now conditioned to disconnect $+$ and $-$ arcs). But one can see from the above proof that negative information in fact makes Lemma \ref{l.area} even more likely. Indeed it makes the event $D_1$ of having a dual crossing in the annulus $A_1$ more likely. 
%The rest is left to to the reader. 
\end{remark}

\begin{remark}
We note that $\log |\EFK{}{\xi}{e^{i\, t\, m}}|$ cannot behave like $- |t|^{\frac {16}{15}}$ as $t\to\infty$ because,
by the Lee-Yang theorem, $\EFK{}{\xi}{e^{i\, t\, m}}$, as a function of complex $t$, has infinitely many zeros, all
purely real. Thus, $\log |\EFK{}{\xi}{e^{i\, t\, m}}|$ must diverge to $-\infty$ at an infinite sequence of real $t$
values (i.e., at the zeros) tending to $\infty$. 
\end{remark}

%\begin{remark}
%In fact by looking at the contribution of the smaller scales $\eps:=2^{-k}$, $k\gg 1$, one finds that these contribute very little
%(it gives something like $(1- (t 2^{-15 k /8})^2)^{2^{2k}}$ which is small, the exponent 2 comes from $\cos(\eta) = 1- \eta^2 \ldots$).
%This should easily lead to 
%\[
%\log |\EFK{}{\xi}{e^{i\, t\, m}}| \asymp - |t|^{\frac {16}{15}}\,,
%\]
%as $t\to\infty$. 
%Well, and yet this is only a heuristics since getting a rigorous lower bound here is non-trivial, it also requires to control what happens at scale $\eps$, so we are far from the lower bound in fact... And even worse: getting an exact constant in this statement seems even much harder here, since ideas like looking at the free energy don't seem to work. 
%\end{remark}

\section{Near-critical magnetization fields}\label{s.NC}

%\section{Massive near-critical magnetization fields}\label{s.NC}
%\section{Near-critical scaling limits}\label{s.NC}

We start by establishing Proposition \ref{pr.NC}.

\ni
{\bf Proof of Proposition \ref{pr.NC}:}

\ni
Let us assume that the boundary condition $\xi$ is $+$ along $\p \Omega$ (the case of free b.c. is treated in the same manner). 
The Ising model with an external field $h_a:= h \, a^{15/8}$ can be thought of as a simple change of measure with respect to the Ising model without external field. In particular, one has for any field $\Phi$:
\begin{align*}
\Pb{\Phi^{a,h}=\Phi} =  \frac{e^{h\, \<{\Phi, 1_\Omega}}}{ \Eb{e^{h \<{\Phi^a, 1_\Omega}}}}  \, \Pb{\Phi^{a,h=0}=\Phi}\,.
\end{align*}
Or, written in terms of the Radon-Nikodym derivative, one has 
\begin{align*}
\frac {d \mu_\Omega^{a,h}} {d\mu_\Omega^a} (\Phi) =  \frac{e^{h\, \<{\Phi, 1_\Omega}}}{ \Eb{e^{h \<{\Phi^a, 1_\Omega}}}}  =  \frac{e^{h\, \<{\Phi, 1_\Omega}}}{ \mu^a_\Omega \bigl[e^{h \<{\Phi, 1_\Omega}}\bigr] } \,,
\end{align*}
where $\mu^{a,h}$ and $\mu^a$ denote respectively the laws of $\Phi^{a,h}_\Omega$ and $\Phi^a_\Omega$.
Now it is not hard to check, using the fact that $\Phi^\infty\sim \mu^{\infty}_\Omega$ has exponential 
moments (Proposition~\ref{pr.expmom}), that $\mu^{a,h}_\Omega$ converges weakly for the topology of $\Sob^{-3}(\Omega)$ to the measure $\mu^{\infty,h}_\Omega$ which is absolutely continuous w.r.t $\mu^\infty_\Omega$ and whose Radon-Nikodym derivative is given by 
\begin{align*}
\frac {d \mu_\Omega^{\infty,h}} {d\mu_\Omega^\infty} (\Phi) =  \frac{e^{h\, \<{\Phi, 1_\Omega}}}{ \Eb{e^{h \<{\Phi^\infty, 1_\Omega}}}} \,.
\end{align*}
We refer to the Appendix A of \cite{CGN1} for details on the topological setup used here ($\Sob^{-3}$).
\qed

We now wish to prove Theorem \ref{th.NC}. It is based on the lemma below together with 
Proposition~\ref{pr.NC}. In what follows, for each $L\in \N$, we will denote 
by $\Lambda_L$  the domain $[-2^L,2^L]^2$.  

\begin{lemma}\label{l.coupling}
For any $\alpha>0$, there exists $L=L(h,\alpha)\in \N$ sufficiently large so that, uniformly in $0<a<\alpha$, one can find a coupling of $\Phi^{a,h}_{\Lambda_L}$ with $\Phi^{a,h}_\C$ satisfying
\begin{align*}
\Eb{ \|\Phi^{a,h}_{\Lambda_L} - \Phi^{a,h}_\C \|_{\Sob^{-3}_\C}} < \alpha\,, 
\end{align*}
where $\| \cdot \|_{\Sob^{-3}_\C}$ is defined by 
\begin{align*}
\| h \|_{\Sob^{-3}_\C}:= \sum_{k\geq 1} \frac 1 {2^k}  \bigl( \| h_{|\Lambda_k} \|_{\Sob^{-3}_{\Lambda_k}}\wedge 1 \bigr)\,.
\end{align*}
\end{lemma}

\begin{remark}
It is easy to check that the distance defined in Lemma \ref{l.coupling} induces the same topology on $\Sob^{-3}_\C$ as the one introduced in Appendix A of \cite{CGN1}. 
\end{remark}

\ni
{\bf Proof of Lemma \ref{l.coupling}:}
Let $L_1\in \N$ be such that $\sum_{k\geq L_1} 2^{-k} < \alpha/2$ (its value will be fixed later, also depending on the value of $h>0$). We wish to find some $L_2 \gg L_1$ such that one can couple the fields $\Phi^{a,h}_{\Lambda_{L_2}}$ and $\Phi^{a,h}_\C$ in such a way that with probability at least $1-\alpha$ they are identical once restricted to the sub-domain $\Lambda_{L_1}$. By the definition of $\|\cdot\|_{\Sob^{-3}_\C}$, this will clearly imply our result with $L=L_2$. 

The coupling will be constructed similarly as in \cite{GPS2a} and \cite{CGN1}. We will also use the FK representation with a ghost vertex used in \cite{CGN2} (see also for example \cite{ghost}). We refer to \cite{CGN2} for more details on this representation. Since the proof below will follow very closely the proof of the lower bound given in Section 3 in \cite{CGN2}, we will not give the full details here.

Following the notation of \cite{CGN2}, let $\bar \omega^{a,h}_\C=(\omega^{a,h}_\C, \tau^{a,h}_\C)$ and $\bar \omega^{a,h}_{L_2}=(\omega^{a,h}_{L_2}, \tau^{a,h}_{L_2})$ be respectively the FK representations of the Ising model with external field $h>0$ on $a \Z^2$ and on $a\Z^2 \cap \Lambda_{L_2}$ with $+$ boundary conditions. These configurations are FK percolation configurations on the graph $a\, \Z^2 \cup \{\g\}$ and the notation $\bar \omega= (\omega,\tau)$ distinguishes between the nearest neighbor edges in $a\, \Z^2$ ($\omega$) and the edges of the type $\<{x,\g}$, with $x\in a\,\Z^2$ ($\tau$).  Furthermore, it is easy to check that $\omega^{a,h}_{L_2}$ stochastically dominates  $\omega^{a,h}_\C$. Let us divide the annulus $A_{L_1,L_2}$ into disjoint annuli of ratio 4: namely $A_1:= A_{4^{-1} L_2, L_2}, A_2:= A_{4^{-2} L_2, 4^{-1} L_2}$ and so on. As such one has about $\log_{4}(L_2/L_1)$ annuli. As in \cite{CGN1}, we will explore ``inward''  both configurations by preserving the monotonicity $\omega_\C \preccurlyeq \omega_{L_2}$ and by trying to find a matching circuit $\gamma$ in each annulus with positive probability. As in \cite{CGN1}, the main ingredient for the coupling is the RSW theorem from \cite{RSWfk}. The difference in our present setting is that one also has to deal with the influence of the {\it ghost} vertex $\g$. In particular, finding a matching circuit $\gamma$ is not enough if one wants to claim that the conditional law ``inside'' 
 $\gamma$ are the same: one also has to make sure that the circuit $\gamma$ is connected in both configurations to the ghost vertex $\g$. 

We proceed as follows. Assume we did not succeed in coupling the two configurations in the first $i-1$ annuli $A_1,\ldots, A_{i-1}$ and consider the $i^{th}$ annulus $A_i = A_{4^{-i} L_2, 4^{-(i-1)L_2}}$. At this point, the configurations have been explored everywhere except inside the outer boundary of $A_i$, %$\p_1 A_i$
and $\omega^{a,h}_{L_2}$ dominates $\omega^{a,h}_{\C}$. Inside the annulus $A_i$, we will distinguish 3 sub-annuli : $B_1 = A_{3. 4^{-i}, 4^{-i+1}}$, $B_2 = A_{2.4^{-i}, 3. 4^{-i}}$ and $B_3=A_{4^{-i}, 2. 4^{-i}}$. From the RSW theorem of \cite{RSWfk}, there are open circuits in $\omega^{a,h}_\C$ (and thus $\omega^{a,h}_{L_2})$ with positive probability $c>0$ in each of $B_1,B_2,B_3$. This is due to the fact that $\omega^{a,h}_\C$ dominates a critical FK configuration with zero magnetic field and with wired boundary conditions along $\p_1 A_i$ (see Section 3 in \cite{CGN2}). Furthermore, due to the positive information inside $B_2$ (thanks to the open circuits in each $B_1$ and $B_3$), it is easy to extend the techniques used to prove Lemma 3.1. in \cite{CGN2} (i.e. an appropriate second moment argument) to show that with positive probability $c>0$, there are at least $c\, (4^{-i} L_2/a)^{15/8}$ points inside $B_2$ which are connected in $\omega^{a,h}_\C$ to the ``outermost'' open circuit $\gamma$ for the configuration $\omega^{a,h}_\C$ in the annulus $B_3$.  Since $(4^{-i} L_2)^{15/8} \geq (L_1)^{15/8}$, the exact same proof as for Lemma 3.2 of \cite{CGN1} shows that if one chooses $L_1$ large enough (depending on $h$), then conditioned on the above event of having at least $c\, (4^{-i} L_2/a)^{15/8}$ points connected to $\gamma$, with conditional probability at least 1/2, the cluster including $\gamma$ will be connected to the ghost vertex $\g$ for the configuration $\omega^{a,h}_\C$ (and thus for $\omega^{a,h}_{L_2}$ as well). Once $\omega^{a,h}_{\C}$ and $\omega^{a,h}_{L_2}$ have a matching circuit $\gamma$ connected to $\g$, one can sample the rest of the configurations so that they match ``inside'' the circuit $\gamma$. (As in \cite{CGN1}, the exploration process is driven by $\omega^{a,h}_\C$.)  To conclude, we choose $L_1$ so that it satsifies the two constraints discussed above (i.e.,
$\sum_{k\geq L_1} 2^{-k} <\alpha/2$ and the constraint relative to $h>0$). This gives a us a certain positive probability $c>0$ to couple both configurations in any annulus $A_{4^{-1}L, L}$, $L\geq L_1$. The proof is then completed by choosing $L=L(h,\alpha)=L_2$ so that $c^{\log_4(L_2/L_1)}<\alpha/2$. \qed

\ni
{\bf Proof of Theorem \ref{th.NC}:}

\ni
By Proposition \ref{pr.NC}, for any $L\in \N$, one has that $\Phi^{a,h}_{\Lambda_L}$ converges in law to $\Phi^{\infty,h}_{\Lambda_L}$ in $\Sob^{-3}_{\Lambda_L}$. It is easy to check that this convergence in law also holds in the space $(\Sob^{-3}_\C, \|\cdot \|_{\Sob^{-3}_\C})$. 
Since this latter space is  Polish, for any $\alpha>0$, there exists $a_0=a_0(\alpha)>0$ such that, for any $a <a_0$, one can couple $\Phi^{a,h}_{\Lambda_L}$ with $\Phi^{\infty,h}_{\Lambda_L}$ so that 
\[
\Eb{ \|\Phi^{a,h}_{\Lambda_L} - \Phi^{\infty,h}_{\Lambda_L} \|_{\Sob^{-3}_\C}} < \alpha\, . 
\]

By using this fact together with Lemma \ref{l.coupling} and the fact that $(\Sob^{-3}_\C, \|\cdot \|_{\Sob^{-3}_\C})$ is Polish, one easily obtains that $\{\Phi^{\infty,h}_{\Lambda_L}\}_{L\in \N}$ converges in law in $\Sob^{-3}_\C$ as $L\to \infty$ to a limiting field $\Phi^{\infty,h}_\C$. Now that our limiting random field is defined, to conclude about the convergence in law of $\Phi^{a,h}_\C$ to this limiting field, we proceed in the same manner: for any $\eps>0$, one can find $a_0>0$ sufficiently small so that for any $a<a_0$, there exists a joint coupling $(\Phi^{a,h}_\C, \Phi^{a,h}_{\Lambda_L}, \Phi^{\infty,h}_{\Lambda_L}, \Phi^{\infty,h}_\C)$ such that all fields are $\eps$-close to each other (for $\|\cdot\|_{\Sob^{-3}_\C}$) with probability at least $1-\eps$. This proves 
the convergence in law of $\Phi^{a,h}_\C$ to $\Phi^{\infty,h}_\C$. 
%\margin{New end of proof below.}

%To complete the proof of the theorem, we need to show
%decay of the two-point correlation for $\Phi^{\infty,h}$. Here is a sketch.
%This is shown uniformly in $a$ for $\Phi^{a,h}$ by again using the FK 
%representation with a ghost vertex.
%There, $\sigma_{x}$ and $\sigma_{y}$ are (conditionally) independent unless $x$ and $y$ are in the
%same FK cluster {\it and\/} that cluster is {\it not\/} connected to the ghost.
%But it is unlikely (exponentially in the distance $||x - y||$)
%that a cluster connects $x$ and $y$ without containing
%$\tilde{c} (||x - y||/a)^{15/8}$ other vertices and then even
%more unlikely for none to connect to the ghost. Estimating the size of a cluster connecting
%$x$ and $y$ uses RSW arguments to connect order $||x - y||$ one-arm vertices near any
%path connecting $x$ and $y$ to that path. 
\qed

%\bibliographystyle{../halpha}
%\bibliographystyle{alpha}
%\addcontentsline{toc}{section}{Bibliography}
%
%%\bibliography{../mr,../prep,../notmr,references}
%\bibliography{magnetization}

\ \\
\ \\
{\bf Federico Camia}\\
Department of Mathematics, Vrije Universiteit Amsterdam and NYU Abu Dhabi. \\
Research supported in part by NWO Vidi grant 639.032.916\\
\\
{\bf Christophe Garban}\\
ENS Lyon, CNRS\\
\url{http://perso.ens-lyon.fr/christophe.garban/}\\
Partially supported by ANR grant MAC2 10-BLAN-0123\\
\\
{\bf Charles M. Newman}\\
Courant institute of Mathematical Sciences, \\
New York University, New York, NY 10012, USA\\
Research supported in part by NSF grants OISE-0730136 and MPS-1007524

\end{document}